\documentclass{article}

\PassOptionsToPackage{numbers, compress}{natbib}

\usepackage[final]{nips_2018}

\usepackage{algorithm}
\usepackage[utf8]{inputenc} 
\usepackage[T1]{fontenc}    
\usepackage{hyperref}       
\usepackage{url}            
\usepackage{booktabs}       
\usepackage{amsfonts}       
\usepackage{nicefrac}       
\usepackage{microtype}      

\usepackage{amsmath,amsfonts,amssymb,amsthm,array}
\usepackage{amsmath}

\usepackage{caption}
\usepackage{subcaption}
\usepackage{algorithm}
\usepackage{algpseudocode}
\usepackage{bm}
\usepackage{color}
\usepackage{graphicx}
\usepackage{cancel}

\usepackage{amsmath,amsfonts,amssymb,amsthm,array}
\usepackage{amsmath}
\usepackage{algorithm}
\usepackage{caption}
\usepackage{subcaption}
\usepackage{algpseudocode,algorithm,algorithmicx}
\usepackage{bm}
\usepackage{color}
\usepackage{graphicx}

\usepackage[algo2e]{algorithm2e} 

\newcommand{\ac}{\alpha}

\newcommand{\E}[1]{{\mathbb{E}\left[#1\right] }}    

\newcommand{\R}{\mathbb{R}}

\usepackage{array}
\newcolumntype{P}[1]{>{\centering\arraybackslash}p{#1}}
\newcolumntype{M}[1]{>{\centering\arraybackslash}m{#1}}

\newcommand{\eqdef}{\overset{\text{def}}{=}}

\newcommand{\cE}{{\cal E}}

\newcommand{\cG}{{\cal G}}

\newcommand{\cV}{{\cal V}}

\newcommand{\bA}{{\bf A}}

\newcommand{\mA}{{\bf A}}

\theoremstyle{plain}
\newtheorem{thm}{Theorem}[]

\newtheorem{defn}[thm]{Definition}
\newtheorem{cor}[thm]{Corollary}

\theoremstyle{remark}

\usepackage[colorinlistoftodos,bordercolor=orange,backgroundcolor=orange!20,linecolor=orange,textsize=scriptsize]{todonotes}

\title{A Privacy Preserving Randomized Gossip Algorithm \\
via Controlled Noise Insertion}

\author{Filip Hanzely$^a$ \And Jakub Kone\v{c}n\'{y}${}^b$  \And Nicolas Loizou${}^{c,}$\thanks{Corresponding author (n.loizou@sms.ed.ac.uk)}
\AND
Peter Richt\'{a}rik${}^{a, c, d}$
\And Dmitry Grishchenko${}^e$\\ 
\phantom{xxx}
\\
\hspace{-160pt}${}^a$ \em King Abdullah University of Science and Technology, KSA\\
\hspace{-160pt}${}^b$ \em Google, USA\\
\hspace{-160pt}${}^c$ \em University of Edinburgh, UK\\
\hspace{-160pt}${}^d$ \em Moscow Institute of Physics and Technology, Russia\\
\hspace{-160pt}${}^e$ \em Universit\'{e} Grenoble Alpes, France
}

\begin{document}

\maketitle

\begin{abstract}
In this work\footnote{The full-length paper, which includes a number of additional algorithms and results (including proofs of statements and experiments), is available in \cite{hanzely2017privacy}.} we present a randomized gossip algorithm for solving the average consensus problem while at the same time protecting the information about the initial private values stored at the nodes. We give iteration complexity bounds for the method and perform extensive numerical experiments. 
\end{abstract}

\section{Introduction}

In this paper we consider the average consensus (AC) problem. Let $\cG=(\cV,\cE)$ be an undirected connected network  with node set $\cV=\{1,2,\dots,n\}$ and edges $\cE$ such that $|\cE| = m$. Each node $i \in \cV$ ``knows'' a private value $c_i \in \R$. The goal of AC is for every node of the network to compute the average of these values, $\bar{c}\eqdef\tfrac{1}{n}\sum_i c_i$, in a distributed fashion. That is, the exchange of information can only occur between connected nodes (neighbours).

The literature on distributed protocols for solving the average consensus problem is vast and has long history \cite{tsitsiklis1984problems, tsitsiklis1986distributed, bertsekas1989parallel, kempe2003gossip}. In this work we focus on one of the most popular class of methods for solving the average consensus, the randomized gossip algorithms and propose a gossip algorithm for protecting the information of the initial values $c_i$, in the case when these may be sensitive. In particular, we develop and analyze a privacy preserving variant of the randomized pairwise gossip algorithm (``randomly pick an edge $(i,j)\in \cE$ and then replace the values stored at vertices $i$ and $j$ by their average'') first proposed in \cite{boyd2006randomized} for solving the average consensus problem. 
While we shall not formalize the notion of privacy preservation in this work, it will be intuitively clear that our methods indeed make it harder for nodes to infer information about the private values of other nodes, {\it which might be useful in practice.} 

\subsection{Related Work on Privacy Preserving Average Consensus} 
The introduction of notions of privacy within the AC problem is relatively recent in the literature, and the existing works consider two different ideas.
In \cite{huang2012differentially}, the concept of differential privacy \cite{dwork2014algorithmic} is used to protect the output value $\bar{c}$ computed by all nodes. In this work, an exponentially decaying Laplacian noise is added to the consensus computation. This notion of privacy refers to protection of the {\it final average}, and formal guarantees are provided. A different line of work with a more stricter goal is the design of privacy-preserving average consensus protocols that guarantee protection of the {\it initial values} $c_i$ of the nodes \cite{nozari2017differentially, manitara2013privacy, mo2017privacy}. In this setting each node should be unable to infer a lot about the initial values $c_i$ of any other node. In the existing works, this is mainly achieved with the clever addition of noise through the iterative procedure that guarantees preservation of privacy and at the same time converges to the exact average. We shall however mention, that none of these works address any specific notion of privacy (no clear measure of privacy is presented) and it is still not clear how the formal concept of differential privacy \cite{dwork2014algorithmic} can be applied in this setting.
\subsection{Main Contributions}
In this work, we present the first randomized gossip algorithm for solving the Average Consensus problem while at the same time protecting the information about the initial values. To the best of our knowledge, this work is the first which combines the \emph{gossip-asynchronous framework} with the privacy concept of protection of the initial values. Note that all the previously mentioned privacy preserving average consensus papers propose protocols which work on the synchronous setting (all nodes update their values simultaneously).

The convergence analsysis of proposed gossip protocol (Algorithm~\ref{PrivacyGossip}) is dual in nature. The dual approach is explained in detail in Section~\ref{sec:duality}. It was first proposed for solving linear systems in \cite{SDA, loizou2017momentum} and then extended to the concept of average consensus problems in \cite{LoizouRichtarik, loizou2018accelerated}. The dual updates immediately correspond to updates of the primal variables, via an affine mapping.

Algorithm~\ref{PrivacyGossip} is inspired by the works of \cite{manitara2013privacy, mo2017privacy}, and protects the initial values by inserting noise in the process. Broadly speaking, in each iteration, each of the sampled nodes first adds a noise to its current value, and an average is computed afterward. Convergence is guaranteed due to the correlation in the noise across iterations. Each node remembers the noise it added last time it was sampled, and in the following iteration, the previously added noise is first subtracted, and a fresh noise of smaller magnitude is added. Empirically, the protection of initial values is provided by first injecting noise into the system, which propagates across the network, but is gradually withdrawn to ensure convergence to the true average.

\section{Technical Preliminaries} 
\label{sec:duality}

\paragraph{Primal and Dual Problems}
Consider solving the (primal) problem of projecting a given vector $c=x^0\in \R^n$ onto the solution space of a  linear system: 
\begin{equation}\label{eq:primal}\min_{x\in \R^n} \{P(x) \eqdef \tfrac{1}{2}\|x-x^0\|^2\} \quad \text{subject to} \quad \bA x=b,\end{equation} 
where $\bA\in \R^{m\times n}$, $b\in \R^m$, $x^0\in \R^n$.
We assume the problem is feasible, i.e., that the system $\mA x = b$ is consistent. With the above optimization problem we associate the dual problem
\begin{equation}\label{eq:dual}\max_{y\in \R^m} D(y)\eqdef (b-\bA x^0)^\top y - \tfrac{1}{2}\|\bA^\top y\|^2.\end{equation}
The dual is an unconstrained concave (but not necessarily strongly concave) quadratic maximization problem. It can be seen that as soon as the system $\mA x = b$ is feasible, the dual problem is bounded. Moreover, all bounded concave quadratics in $\R^m$ can be written in the as $D(y)$ for some matrix $\mA$ and vectors $b$ and $x^0$ (up to an additive constant).

With any dual vector $y$ we associate the primal vector via an affine transformation, $x(y) = x^0 + \mA^\top y.$
It can be shown that if $y^*$ is dual optimal, then $x^*=x(y^*)$ is primal optimal. Hence, any dual algorithm producing a sequence of dual variables $y^t \to y^*$ gives rise to a corresponding primal algorithm producing the sequence $x^t \eqdef x(y^t) \to x^*$. See \cite{SDA, loizou2017momentum} for the correspondence between primal and dual methods.

\paragraph{Randomized Gossip Setup: Choosing $\mA$.}
In the gossip framework we wish $(\mA,b)$ to be an {\em average consensus (AC)} system.
\begin{defn} (\cite{LoizouRichtarik}) 
Let  $\cG = (\cV, \cE)$ be an undirected graph with $|\cV| = n$ and $|\cE|=m$. Let $\bA$ be a real matrix with $n$ columns.
The linear system $\bA x = b$ is an ``average consensus (AC) system'' for graph $\cG$ if $\bA x = b$ iff $x_i = x_j$ for all $(i,j) \in \cE$.
\end{defn}
In the rest of this paper we focus on a specific AC system; one  in which the matrix $\bA$  is the incidence matrix of the graph $\cG$ (see Model 1 in \cite{SDA}).  In particular, we let $\mA\in \R^{m\times n}$ be the matrix defined as follows. Row $e=(i,j)\in \cE$  of $\mA$ is given by $\mA_{ei} = 1$, $\mA_{ej}=-1$ and $\mA_{el}=0$ if $l\notin \{i,j\}$.  Notice that the system $
\mA x=0$ encodes the constraints $x_i=x_j$ for all $(i,j)\in \cE$, as desired. It is also known that randomized Kaczmarz method~\cite{RK, gower2015randomized,loizou2017linearly} applied to Problem~\ref{eq:primal} is equivalent to randomized gossip algorithm (see \cite{LoizouRichtarik, loizou2018accelerated, loizou2018provably} for more details). 

\section{Private Gossip via Controlled Noise Insertion}
 \label{sec:Private}

In this section, we present the Gossip algorithm with Controlled Noise Insertion. As mentioned in the introduction, the approach is similar to the technique proposed in \cite{manitara2013privacy, mo2017privacy}. Those works, however, address only algorithms in the synchronous setting, while our work is the first to use this idea in the asynchronous setting. Unlike the above, we provide finite time convergence guarantees and allow each node to add the noise differently, which yields a stronger result.

In our approach, each node adds noise to the computation independently of all other nodes. However, the noise added is correlated between iterations for each node. We assume that every node owns two parameters --- the initial magnitude of the generated noise $\sigma_i^2$ and rate of decay of the noise $\phi_i$. The node inserts noise $w_i^{t_i}$ to the system every time that an edge corresponding to the node was chosen, where variable $t_i$ carries an information how many times the noise was added to the system in the past by node $i$. Thus, if we denote by $t$ the current number of iterations, we have $\sum_{i=1}^nt_i = 2t$.

In order to ensure convergence to the optimal solution, we need to choose a specific structure of the noise in order to guarantee the mean of the values $x_i$ converges to the initial mean. In particular, in each iteration a node $i$ is selected, we subtract the noise that was added last time, and add a fresh noise with smaller magnitude:
$w_i^{t_i} = \phi_i^{t_i}v_i^{t_i}-\phi_i^{t_i-1}v_i^{t_i-1},$ where $0 \leq \phi_i<1, v_i^{-1} = 0$ and $v_i^{t_i}\sim N(0,\sigma_i^2)$ for all iteration counters $k_i \geq 0$ is independent to all other randomness in the algorithm. This ensures that all noise added initially is gradually withdrawn from the whole network.

After the addition of noise, a standard Gossip update is made, which sets the values of sampled nodes to their average. Hence, we have $\lim_{t\rightarrow \infty}\E{\left(\overline{c}-\frac1n\sum_{i=1}^n x_i^t\right)^2}= 0,$ as desired. 

It is not the purpose of this paper to define any quantifiable notion of protection of the initial values formally. However, we note that it is likely the case that the protection of private value $c_i$ will be stronger for bigger $\sigma_i$ and for $\phi_i$ closer to $1$.

\begin{algorithm}[t!]
\caption{Privacy Preserving Gossip Algorithm via Controlled Noise Insertion}
\label{PrivacyGossip}
\textbf{Input: }{vector of private values $c\in \R^n$; initial variances $\sigma^2_i \in \R_+$ and variance decrease rate $\phi_i$ such that $0\leq \phi_i < 1$ for all nodes $i$.}\\
{\textbf{Initialize}:} Set $x^0=c$; $t_1=t_2=\dots = t_n=0$, $v_1^{-1}=v_2^{-1}=\dots = v_n^{-1}=0$.\\
\For {$t= 0,1,\dots k-1$} {
     \begin{enumerate}
\item Choose edge $e = (i,j)\in \cE$ uniformly at random
\item Generate $v_i^{t_i}\sim N(0,\sigma^2_i)$ and $v_j^{t_j}\sim N(0,\sigma^2_j)$
\item Set $w_i^{t_i} = 
 \phi_i^{t_i}v_i^{t_i}-\phi_i^{t_i-1}v_i^{t_i-1} $
  and 
 $w_j^{t_j} = 
 \phi_j^{t_j}v_j^{t_j}-\phi_j^{t_j-1}v_j^{t_j-1} 
 $
\item Update the primal variable: $x^{t+1}_i =x^{t+1}_j= 
\tfrac{x^t_i+w^{t_i}_i+x^t_j+w^{t_j}_j}{2}$, $\forall\, l \neq i,j:\,x^{t+1}_l=x^{t}_l
  $
\item Set $t_i=t_i+1$ and $t_j=t_j+1$
\end{enumerate}
 }
\textbf{return} $x^k$
\end{algorithm}

We now provide results of dual analysis of Algorithm~\ref{PrivacyGossip}.

\begin{thm}
Let us define $\rho \eqdef 1-\frac{\ac(\cG)}{2m}$ and $
\psi^t \eqdef \frac{1}{\sum_{i=1}^n\left(d_i\sigma_i^2\right)}\sum_{i=1}^n d_i \sigma_i^2\left(1-\frac{d_i}{m}\left(1-\phi_i^2\right) \right)^{t},$ where $\ac(\cG)$ stands for algebraic connectivity of $\cG$ and $d_i$ denotes the degree of node $i$.
Then for all $k\geq 1$ we have the following bound
\begin{equation*}
\E{ D(y^*)- D(y^{k}) } \leq  \rho^k \left( D(y^*)- D(y^{0}) \right)  + \frac{\sum\left(d_i\sigma_i^2\right)}{4m}\sum_{t=1}^k \rho^{k-t}\psi^{t}.
\end{equation*}
\label{T: ng general convergence}
\end{thm}

Note that $\psi^t$ is a weighted sum of $t$-th powers of real numbers smaller than one.  For large enough $t$, this quantity will depend  on the largest of these numbers. This brings us to define $M$ as the set of indices $i$ for which the quantity $1-\frac{d_i}{m}\left(1-\phi_i^2\right)$ is maximized: $M=\arg\max_{i} \left\{ 1-\frac{d_i}{m}\left(1-\phi_i^2\right)\right\}. $
Then for any $i_\mathrm{max}\in M$ we have
$$
\psi^t  
\approx
 \frac{1}{\sum_{i=1}^n \left(d_i\sigma_i^2\right)} \sum_{i\in M} d_i \sigma_i^2\left(1-\frac{d_i}{m}\left(1-\phi_i^2\right) \right)^{t}
 =
  \frac{\sum_{i\in M} d_i \sigma_i^2}{\sum_{i=1}^n \left(d_i\sigma_i^2\right)} \left(1-\frac{d_{i_\mathrm{max}}}{m}\left(1-\phi_{i_\mathrm{max}}^2\right) \right)^{t},
$$
which means that increasing $\phi_j$ for $j\not\in M$ will not substantially influence convergence rate. Note that as soon as we have
\begin{equation}
\label{eq: treshold}
\rho> 1-\frac{d_{i}}{m}\left(1-\phi_{i}^2\right)  
\end{equation} 
for all $i$, the rate from theorem \ref{T: ng general convergence} will be driven by $\rho^k$ (as $k \rightarrow \infty$) and we will have $ \E{ D(y^*)- D(y^{k}) } = \tilde{O}\left(\rho^k\right)$.
One can think of the above as a threshold: if there is $i$ such that $\phi_i$ is large enough so that the inequality \eqref{eq: treshold} does not hold, the convergence rate is driven by $\phi_{i_\mathrm{max}}$. Otherwise, the rate  is not influenced by the insertion of noise. Thus, in theory, we do not pay anything in terms of performance as long as we do not hit the threshold. One might be interested in choosing $\phi_i$ so that the threshold is attained for all $i$, and thus $M=\{1,\dots,n\}$. This motivates the following result:
\begin{cor}
\label{corolary}
Let us choose $\phi_i \eqdef \sqrt{1-\frac{\gamma}{d_i}}$ for all $i$, where  $\gamma \leq d_{\mathrm{min}}$.
Then
\begin{eqnarray*}
\E{ D(y^*)- D(y^{k}) } & \leq &
\left( 1-\min\left( \frac{\ac(\cG)}{2m},\frac{\gamma}{m}\right) \right)^k \left( D(y^*)- D(y^{0})+\frac{\sum_{i=1}^n \left(d_i\sigma_i^2\right)}{4m}k \right). 
\end{eqnarray*}
As a consequence, $\phi_i=\sqrt{1-\tfrac{\ac(\cG)}{2d_i}}$ is the largest decrease rate of noise for node $i$ such that the guaranteed convergence rate of the algorithm is not violated.
\label{C: noisy gossip special}
\end{cor}

\section{Experiments}
\label{exper}
In this section we present a preliminary experiment (for more experiments see Section~\ref{moreExp}, in the Appendix) to evaluate the performance of the Algorithm~\ref{PrivacyGossip} for solving the Average Consensus problem. The algorithm has two different parameters for each node $i$. These are the initial variance $\sigma_i^2 \geq 0$ and the rate of decay, $\phi_i$,  of the noise.

In this experiment we use two popular graph topologies the cycle graph (ring network) with $n=10$ nodes and the random geometric graph with $n=100$ nodes and radius $r=\sqrt{\log(n)/n}$.

In particular, we run Algorithm~\ref{PrivacyGossip} with $\sigma_i = 1$ for all $i$, and set $\phi_i = \phi$ for all $i$ and some $\phi$. We study the effect of varying the value of $\phi$ on the convergence of the algorithm. 

\begin{figure}[H]
\centering
\begin{subfigure}{.45\textwidth}
  \centering
  \includegraphics[width=1\linewidth]{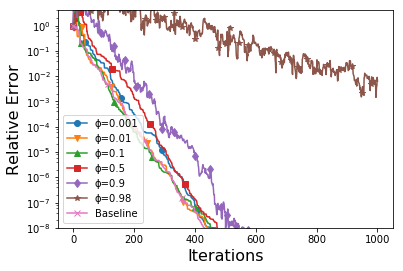}
  \caption{}
\end{subfigure}
\begin{subfigure}{.45\textwidth}
  \centering
  \includegraphics[width=1\linewidth]{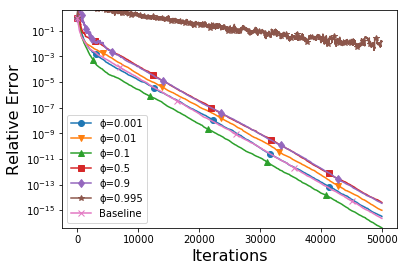}
  \caption{}
\end{subfigure}
\caption{Convergence of Algorithm~\ref{PrivacyGossip}, on the cycle graph (left) and random geometric graph (right) for different values of $\phi$. The ``Relative Error " on the vertical axis represents the $\frac{\|x^t-x^*\|^2}{\|x^0 - x^*\|^2}$}
\label{plot}
\end{figure}

In Figure~\ref{plot} we see that for small values of $\phi$, we eventually recover the same rate of linear convergence as the Standard Pairwise Gossip algorithm  (Baseline) of \cite{boyd2006randomized}. If the value of $\phi$ is sufficiently close to $1$ however, the rate is driven by the noise and not by the convergence of the Standard Gossip algorithm. This value is $\phi = 0.98$ for cycle graph, and $\phi=0.995$ for the random geometric graph in the plots we present.

\small
\bibliographystyle{plain}
{\footnotesize\bibliography{PPMLarxiv}}
\end{document}